\documentclass[onecolumn,aps,prb,showpacs,
 amsfonts,amssymb,amsmath]{revtex4}
\usepackage{mathrsfs}
\usepackage{amsmath}
\usepackage{graphicx}

\begin{document}
\draft
\title{Regularization to orthogonal-polynomials fitting with application to magnetization data}

\author{Bao Xu$^{1,2}$}

\affiliation{$^{1}$ Department of Physics Science and Technology, 
Baotou Teachers' College, Baotou 014030, China}

\affiliation{$^{2}$ Key Laboratory of Magnetism and Magnetic Materials 
 at Universities of Inner Mongolia Autonomous Region, Baotou 014030, China}
 

\footnotetext{\hspace*{-.45cm}\footnotesize Project supported by National
Natural Science Foundation of China (Grant No. 61565013) and The Scientific Research Project of the Inner Mongolia Autonomous Region Colleges and Universities (Grant No. NJZZ199).}
\footnotetext{\hspace*{-.45cm}\footnotesize  Corresponding author. E-mail: xubao79@163.com }

\begin{abstract}
An obstacle encountered in applying orthogonal-polynomials fitting 
is how to select out the proper fitting expression.
By adding a Laplace term to the error expression and 
introducing the concept of overfitting degree,
a regularization and corresponding cross validation scheme 
is proposed for two-variable polynomials fitting.
While the Fortran implementation of above scheme is applied to magnetization data, 
a satisfactory fitting precision is reached, 
and overfitting problem can be quantitatively assessed, 
which therefore offers the quite reliable base 
for future comprehensive investigations of magnetocaloric and 
phase-transition properties of magnetic functional materials.
\end{abstract}


\maketitle

%
%
\section{Introduction}  

Surface fit using two-variable orthogonal polynomials has reasonable advantages.$^{[1]}$ 
One among many is not necessary to solve a normal equation which is probably 
ill-conditioned. A second one usually mentioned in literatures is that,
the coefficients of orthogonal polynomials in the fitting expression
do not depend on those of preceding polynomials.
Hence, fit with two-variable orthogonal polynomials has been widely used in the field such as image processing. 
While applying it to magnetization data, 
another advantage that might be worth mentioning but rarely mentioned 
is rooted in physics. When pressure in the sample chamber does not change,  
magnetization of the sample, $M$, as one of the thermodynamic functions,  
can be completely determined as a function of 
two independent parameters, magnetic field $H$ and temperature $T$.  

Application of this method to magnetization data is not new.
For example, a two-step category has been employed 
to fit the magnetization data of gadolinium, and the fitting expression 
is subsequently used to estimate the corresponding magnetic entropy change.$^{[2]}$
The main features of our method are as follow.

Firstly, in recursively generating orthogonal polynomials 
with classical or modified Gram-Schmidt schemes,
orthogonality of the polynomials will progressively deteriorate.
So an iterating scheme is selected in order to preserve the orthogonality.
(Generally speaking, as has been pointed out before, 
one more orthogonalization process is sufficient to 
significantly improve the orthogonality.)

Secondly, a regularization method is introduced. 
In previous fittings to magnetization data, the overfitting problem 
has rarely been taken into account.
The reason for this probably originates in that not very-high-order polynomials 
are included in the final fitting expression.    
However, our numerical results suggest that fluctuations 
indeed appear between experimentally recorded data, 
and such fluctuations become even severer near to the boundaries. 
So, we are convinced that overfitting occurs and has to be addressed. 
That's why the regularization scheme is introduced to relieve the probable overfitting. 
Although the specific formulations are different,  
the present idea of regularization is similar to that used in reference,$^{[3]}$ 
which uses gridding and interpolation to obtain the normal equation to be solved.

Thirdly, a cross-validation scheme is employed to select out the proper regularization parameter.
The main idea is to divide the experimentally recorded data into three groups, namely, the training, cross-validation, and test groups. Among them, the training group includes the most elements, 
and the other two have less data. For each of different regularization parameters,
use the training group to generate orthogonal polynomials and determine 
the unknowns in the fitting expression by minimizing training error;
apply the cross-validation group to selecting out the optimal regularization parameter 
that corresponds to the minimum of cross-validation error;     
and assess the applicability of the determined expression by computing test error.
For orthogonal polynomials fitting, it's found that
the common cross-validation method can not efficiently select 
out the proper regularization parameter. To overcome this problem, 
the concept of overfitting degree is introduced and used to monitor fitting performance.    
 
Fourthly, two methods to implement sampling are provided.
In principle, a good fitting expression should not be affected by the sample method 
to obtain the training group.
However, our results show an apparent difference 
in the cross-validation performance, 
which is associated with the sampling with respect to 
magnetic field or temperature.       
The reason for the different performance 
is attributed to the (dimensionless and size-normalized)  averaged increment 
of recorded magnetic fields unequal to that of recorded temperatures
within the experimental measurement range.

Fifthly, an extended-precision version of the algorithm is used.
On completely determining the fitting expression,
the magnetization at particular magnetic field and temperature 
calculated from the orthogonal polynomials,
should not be different from the value estimated 
from the linearly independent functions, which  
are used to generate those orthogonal polynomials.  
However, in performing a double-precision version of the algorithm,
we noticed in practice that the values 
obtained by the two methods are actually different. 
So, during recursively computing the values 
of linearly independent functions and subsequently generating orthogonal polynomials,
the error accumulation is significant.
That's why we implement the algorithm at an extended-precision level 
in spite of more computational time.    

Sixthly, an equally-spaced distribution of both $H_{i}$ and $T_{j}$
is not a prerequisite to fit the magnetization data $M(H_{i},T_{j})$. 
In principle, our method is also suitable for 
magnetization data randomly distributed over the $H$-$T$ plane.

Finally, after obtaining the fitting expression, 
it's ready to estimate the magnetocaloric quantities 
at arbitrary magnetic field and temperature within the measuring range, and 
comprehensively investigate magnetic-phase-transition properties
including the order of transition, critical exponents, anomalous specific heat, and so forth.    

The rest of this paper is organized as follows:
Next section provides the general formulations and key algorithms used in this work;
Section 3 applies the algorithm to magnetization data; 
and conclusions are put in Section 4.

\section{General formulation and key algorithms}  

Following formulations are quite general 
and not restricted to magnetization data. 
So the experimental data are denoted by 
as $x_i$, $y_i$ and $z_i$ 
instead of $H_{i}$, $T_{i}$ and $M_{i}$.  

\subsection{Two-variable orthogonal polynomials}  

For experimental data $(x_{i},y_{i},z_{i})$ with $i=1,2,\cdots,N$, 
two-variable orthogonal polynomials ${P}_{s}(x,y)$ 
are defined as satisfy

\begin{equation}
\sum_{i=1}^{N} {P}_{t}(x_i,y_i) {P}_{s}(x_i,y_i) = \delta_{t,s} \sum_{i=1}^{N} 
\left( {P}_{t}(x_i,y_i) \right)^2. \label{e1}
\end{equation}
where, $\delta_{t,s}$ is the Kronecker $\delta$-symbol.  
By using the Gram-Schmidt orthogonalization process 
(classical or modified scheme)
to the linearly independent functions
$$ 
h_{t}(x,y) = 1,
x,
y,
x^{2},
xy,
y^{2},
\cdots
$$
with integer numbers $0\leq t\leq L$.
${P}_{s}(x,y)$ can be recursively generated as follows
\begin{equation*}
{P}_{s}(x,y)= a_{ss} h_{s}(x,y) + \sum_{t=0}^{s-1}a_{st} {P}_{t}(x,y),~0\leq s\leq S.
\end{equation*}
In above expression, $S$ is the largest index of the generated orthogonal polynomials.
It is smaller than or equal to the number of created linear independent function $L$.   
If we assign $a_{s,0}=1$, coefficients $a_{ss}$ 
and $a_{st}$ $( 1 \leq t < s \leq S)$ are determined as
\begin{equation}
a_{ss} = - \frac{
\sum_{i=1}^{N} \left(  {P}_{0}(x_{i},y_{i})  \right)^{2}
}{
\sum_{i=1}^{N} {P}_{0}(x_{i},y_{i})h_{s}(x_{i},y_{i}) 
}, \label{e2}
\end{equation}
\begin{equation}
a_{st} =
- a_{ss}\frac{
\sum_{i=1}^{N} {P}_{t}(x_{i},y_{i}) h_{s}(x_{i},y_{i})
}{
\sum_{i=1}^{N} \left(  {P}_{t}(x_{i},y_{i}) \right)^{2}
}. \label{e3}
\end{equation}

Note that the orthogonal polynomials defined above are not normalized.
The normalized version can be conveniently obtained through dividing coefficients 
$a_{st}$ $(t=0,\cdots,s)$ by $2$-norm
$\left[ \sum_{i=1}^{N} \left( {P}_{t}(x_{i},y_{i}) \right)^{2} \right]^{1/2}$ as
\begin{equation}
a_{ss} = - \left[  
\sum_{i=1}^{N} {P}_{0}(x_{i},y_{i})h_{s}(x_{i},y_{i}) \right]^{-1}, \label{e3}
\end{equation}
\begin{equation}
a_{st} = - a_{ss} \sum_{i=1}^{N} {P}_{t}(x_{i},y_{i}) h_{s}(x_{i},y_{i}). \label{e4}
\end{equation}

\subsection{Coefficients of orthogonal polynomials in the fitting expression} 

After the orthogonal polynomials are obtained,
we can expand the fitting expression $f(x,y)$ with $P_{s}(x,y)$ as  

\begin{equation}
f(x,y) = \sum_{t=0}^{S} b_{t}P_{t}(x,y),\label{fit_expr}
\end{equation}
where $S\leq L$ denotes the maximum index of orthogonal polynomials 
used in the fitting expression. 
By minimizing fitting error 
\begin{equation}
\sigma_{1} = \frac{1}{N}
\sum_{i=1}^{N}\left[  f(x_{i},y_{i})-z_{i} \right]^2 \label{sigma1}
\end{equation}
the coefficient of normalized orthogonal polynomials in the fitting expression is determined as
\begin{equation}
b_{t} = \sum_{i=1}^{N} z_{i}P_{t}(x_{i},y_{i}) ,
~~~(0\leq t \leq S). \label{bt1}
\end{equation}

\subsection{Iterating orthogonalization}

Using above orthogonalization processes, 
it's found that orthogonality becomes poorer and poorer. 
Although it's better than the the classical (CGS) scheme,
performance of the modified Gram-Schmidt (MGS) scheme unavoidably becomes poor  
with increasing the largest index of orthogonal polynomials, $S$.
Since the orthogonality is closely related to fitting precision,
we use the following iterating scheme (IGS) to improve the orthogonality. 


\textbf{Step 1} Recursively compute the values of linearly independent functions 
$h(i,s)$ ($i=1,...,N$; $s=0,...,L$).

(Assume that the first $s$ polynomials $P_{t}(:)$ with $0 \leq t \leq s-1$ 
have been orthonormalized and assigned to $h(:,t)$ with $0 \leq t \leq s-1$.
Estimate the value of the $s+1$-th orthonormalized polynomial $P_{s}(:)$ and assign 
it to $h(:,s)$; and save coefficients $a(s,t)$ with $0 \leq t \leq s$ and $b(s)$. 
)


\textbf{Step 2}  Re-orthogonalize $h(:,s)$ and update $a(:,:)$.

(2-a) Compute the modification coefficient
$\delta(t) = \langle h(:,s);h(:,t)\rangle$, $0 \leq t\leq s-1$;

(2-b) Re-orthogonalize
$h(:,s) = h(:,s)-\sum_{t=0}^{s-1}\delta(t)\cdot h(:,t)$;

(2-c) Update coefficients
$a(s,t) = a(s,t)+\delta(t)$, $0\leq t\leq s-1$;

(2-d) Judge whether the orthogonality criterion is satisfied.
If true then continue; else go back to (2-a).   


\textbf{Step 3} Normalize $h(:,s)$ and update $a(:,:)$.

(3-a) Compute the 2-norm
$p=||h(:,s)||_{\rm 2}$;

(3-b) Normalize $h(:,s)$ as 
$h(:,s) = h(:,s)/p$;

(3-c) Update coefficients
$a(s,s) = p$;

(3-d) Update coefficients
$a(s,t) = - a(s,t)\cdot a(s,s)$, $0\leq t\leq s-1$.


\textbf{Step 4} Compute coefficients of $h(:,s)$ in the fitting expression 
$b(s) = \langle h(:,s);f(:)\rangle$.


\textbf{Step 5} Judge whether fitting precision matches the criterion.
If true then $S = s$ and break out the loop; else continue.

\textbf{Step 6} Subtract the projection of $h(:,t)$ from $h(:,s)$ and update $a(:,:)$.


(6-a) Compute the coefficients of subsequent orthogonal polynomials
$a(t,s) = \langle h(:,t) ; h(:,s) \rangle$, $s+1\leq t \leq L$;

(6-b) Subtract the projection of $h(:,t)$ from $h(:,s)$
$h(:,t) = h(:,t) -  a(t,s) \cdot h(:,s)$, $s+1\leq t \leq L$.


\textbf{Step 7} Update $s$ and go back to \textbf{Step 2}.


\subsection{Overfitting and regularization} 

If the largest index $S$ is not a big number, 
then the changing tendency of experimental data 
can not be properly reflected. 
Now fitting error is large and underfitting happens. 
For decreasing fitting error, more orthogonal polynomials are 
successively generated and added to the fitting expression 
until the wished precision is achieved.  
(The more orthogonal polynomials, the higher fitting precision.)
However, too many polynomials will lead to strong local fluctuations in the fitted surface,
and overfitting happens. The reason for this is that higher-order polynomials 
generally imply more inflection points. The degree of overfitting can be controlled by 
regularization like adding so-called penalty functions into the error expression, 
in order to strengthen the stiffness of fitted surfaces.

In contrast to the method of using penalty functions, 
we implement the regularization by adding a Laplace term to the error expression.
In essence, the Laplace method aims at suppressing the changing rate of curve slope.
After adding a Laplace term 

\begin{equation*}
\nabla^{2}f(x,y) = 0, \label{laplace1}
\end{equation*}
with regularization parameter $\lambda$,
the error expression (\ref{sigma1}) is rewritten as

\begin{equation}
\sigma_{2} = \frac{1}{N}\sum_{i=1}^{N}
\left[
f(x_{i},y_{i})-z_{i}
\right]^2
+ \lambda 
\left[
\nabla^{2} f(x_{i},y_{i})
\right]^2.\label{sigma2}
\end{equation}

It's easy to see that the Laplace term in (\ref{sigma2}) affects only coefficients $b_{t}$.
By minimizing (\ref{sigma2}), it is obtained that    

\begin{equation}
b_{t} =
\frac{
- \lambda R_{t} Q_{t} + \sum_{i=1}^{N}z_{i} P_{t}(x_{i},y_{i})  
}{
\lambda  \left( Q_{t} \right)^{2} + 1 },			\label{bt2}
\end{equation} 
where,

\begin{eqnarray*}
R_{t} &=& \sum_{r=0}^{t-1}b_{r}Q_{r},\\
Q_{t} &=& \sum_{i=1}^{N}\nabla^{2} P_{t}(x_{i},y_{i})
\end{eqnarray*}
with $0\leq t\leq S$.
We next examine whether the Laplace method above really leads to regularization.
Firstly, when $\lambda=0$, equation(\ref{bt2}) reduces to the non-regularized case (\ref{bt1}).
Secondly, if $t\leq 2$, the Laplace term has no contribution to  $b_{t}$ 
since $Q_{0} = Q_{1} = Q_{2} = 0$  (corresponding to linear fitting).
Thirdly, the Laplace term starts to play its role when $t\geq 3$.  
If $\lambda$ is large enough so that equations (\ref{a1}) and (\ref{a2}) are satisfied, 

\begin{eqnarray}
&& \lambda |R_{t}Q_{t}| \gg |\sum_{i=1}^{N}z_{i}P_{t}(x_{i},y_{i})|, \label{a1}
\\ 
&& \lambda  \left( Q_{t} \right)^{2} \gg \sum_{i} \left(P_{t}(x_{i},y_{i})\right)^{2}; 
\label{a2}
\end{eqnarray}
then $b_{t}$ reduces to

\begin{equation*}
b_{t} 
\approx - \frac{\lambda R_{t} Q_{t}}{\lambda (Q_{t})^{2}} 
= - \frac{R_{t}}{Q_{t}}
= - \frac{ \sum_{r=0}^{t-1}b_{r}Q_{r} }{ Q_{t} },
\end{equation*}
namely,
\begin{equation*}
R_{t+1} = \sum_{r=0}^{t} b_{r} Q_{r} \approx 0.
\end{equation*}

If we assume that $R_{t} \approx 0$ when $t\geq t_{0}$,
then $R_{t+1} = b_{t} Q_{t} + R_{t} \approx b_{t} Q_{t} \approx 0$,
which implies that $b_{t}\approx 0$ since $Q_{t}\neq 0$.
Hence, the Laplace term introduced above 
makes $b_{t}$ rapidly decay with increasing $t$, 
so that overfitting is avoided and regularization is realized.

\subsection{Cross validation and overfitting degree} 

Now, we invoke a cross-validation process 
to select out a proper regularization parameter $\lambda$. 
Dividing the whole data into three groups,
one of which includes much more data, labelled ``training group'',
and the other two has fewer data, 
labelled ``cross-validation group" and ``test group", respectively.
For each fixed value of $\lambda$, 
the training group is used to determine coefficient $b_{s}$ by minimizing 
$\sigma_{2}$ in (\ref{sigma2}), and the corresponding training error is computed as

\begin{equation}
\sigma_{\rm tr} = \frac{1}{ N_{\rm train} }
\sum_{i=1}^{ N_{\rm train} }\left[
f(x_{i},y_{i})-z_{i}\right]^2. \label{err_tr}
\end{equation}
Subsequently, the cross-validation group selects out the value of $\lambda$ that minimizes 
the cross-validation error 

\begin{equation}
\sigma_{\rm cv} = \frac{1}{N_{\rm cv}}
\sum_{j=1}^{N_{\rm cv}}\left[
f(x_{i},y_{i})-z_{i}\right]^2. \label{err_cv}
\end{equation}
Finally, the test group assesses the applicability 
of the determined fitting expression by
calculating the test error
\begin{equation*}
\sigma_{\rm test} = \frac{1}{N_{\rm test}}
\sum_{j=1}^{N_{\rm test}}\left[
f(x_{i},y_{i})-z_{i}\right]^2.
\end{equation*} 

In contrast to the ordinary scheme used in the field such as machine learning,
the model used here has two parameters that require determining, 
namely, the number of orthogonal polynomials $S$ and 
the regularization parameter $\lambda$.

By setting certain routine-terminating criterion, the optimal choice of $S$ can be determined by
minimizing the training error $\sigma_{\rm tr}$. 
A useful criterion can be defined by assessing the changing tendency of fitting error.
For example, on increasing $S$ from $S_{0}$, 
if the fitting error is not apparently decreased, 
it is reasonable to consider $S_{0}$ to be the optimal value of $S$.       

Another task is to determine parameter $\lambda$.
Practically, we find that both $\sigma_{\rm tr}$ and $\sigma_{\rm cv}$ 
decrease with enlarging S at fixed $\lambda$ or with increasing $\lambda$ at fixed $S$. 
Thus, the optimal value of $\lambda$ can not be identified to the one 
that minimizes $\sigma_{\rm cv}$.

This motivates us to construct a new quantity  
to characterize the degree of overfitting (and also underfitting).        
Typically, when $\sigma_{\rm cv}$ is approximately equal to $\sigma_{\rm tr}$, 
underfitting happens; if $\sigma_{\rm cv}$ is much larger than $\sigma_{\rm tr}$, 
overfitting occurs. We can define overfitting degree $\gamma$ as
\begin{equation}
\gamma = \ln \Big | 
\frac{ \sigma_{\rm cv}- \sigma_{\rm tr} }{ \sigma_{\rm tr} }
\Big |. \label{overfit1}  
\end{equation}
It is identified as underfitting if $\gamma \ll -1$, and 
overfitting when $\gamma \gg 1$. 

\subsection{Size normalization}

It's noticed that the calculated value of fitting error depends 
on the measurement unit used for experimentally recorded data 
$(X_{i},Y_{i},Z_{i})$ $(i=1,\cdots,N)$ where
$X_{i} \in [X_{\rm min},X_{\rm max}]$,
$Y_{i} \in [Y_{\rm min},Y_{\rm max}]$,
and $Z_{i} \in [Z_{\rm min},Z_{\rm max}]$. 
For comparison purposes, we do size normalization as follow
\begin{eqnarray*}
X &=& f_{X}(x) = X_{\rm min} + x (X_{\rm max} - X_{\rm min}),
\\
Y &=& f_{Y}(y) = Y_{\rm min} + y (Y_{\rm max} - Y_{\rm min}),
\\
Z &=& f_{Z}(z) = Z_{\rm min} + z (Z_{\rm max} - Z_{\rm min}). 
\end{eqnarray*}
where, $x,y,z \in [0,1]$.
Hence, the fitting error is calculated from the after-transformed data. 
With an inverse transformation, 
the physical quantities are obtained in the measurement unit.

\subsection{Uniform sampling}

In dividing the experimental data into three groups,
a uniform sampling algorithm is executed in order to optimize fitting performance. 
Original experimental data are firstly sorted in terms of 
a sample parameter, $X$ or $Y$.
Then the uniform sample is executed on a pro-rata basis,
which is regulated by the sampling factor, and
the sorted data are put into the training, cross-validation, 
and test groups according to the sampling factor.
Numerical results show that difference appears 
between different sampling methods, which is 
attributed to the different data density along 
$X$ axis with that along $Y$ axis.

\subsection{Coefficients of linear independent functions $h_{i}(x,y)$ in the fitting expression} 

After coefficients $a_{st}$ and $b_{s}$ being determined, 
the functional value at arbitrary location $(x,y)$ within the measuring range 
can be estimated in a similar way to that used in the fitting routine.
Another method is to estimate from linearly independent functions as
\begin{equation}
f(x,y) = \sum_{t=0}^{S} c_{t}h_{t}(x,y),  	\label{fit_expr_2}
\end{equation}
where $c_{t}$ is the coefficient of the $t$-th linearly independent functions 
in the fitting expression, and can be readily calculated from $a_{ts}$ and $b_{t}$.
The latter scheme is recommended for lower computation cost.

Practically, in implementing a double-precision version of the algorithm,
it's found that the fitting value calculated from linearly independent functions 
significantly differs from that computed from orthogonal polynomials, 
when the power exponent of the fitting expression is very high.
After an extended-precision algorithm is applied, 
the difference decreases.  
If the extended-precision operation is also used to 
recursively generate linearly independent functions,
the difference is not longer obvious.    
These facts suggest that round-off error is rapidly accumulated while 
recursively computing the linearly independent functions. 
Hence, it is required to carefully consider the influence of error accumulations 
in fitting using two-variable orthogonal polynomials.

\subsection{Some useful recursive formula}

Since the fitting expression is essentially the linear combination 
linearly independent functions $h(x,y)$,
utilizing recursive properties of the corresponding partial derivative and integral,  
it's quite convenient to compute the physical quantities of interest. 
Followings are typical recursive algorithms at given $x$ and $y$ used in this work.    

\subsubsection{$h(i)\leftarrow h_{i}(x,y)$ used to compute linear independent functions as well as magnetization}

\begin{eqnarray}
&&h(0) = 1;~~ h(1) = x;~~ h(2) = y;~~s=1;
\nonumber\\
&&\mathrm{For}~~ m\geq 2
\nonumber\\
&& \{ s = s + m;
\nonumber\\
&&h(s+j) = \left\{
\begin{array}{cl}
x\cdot h(s-m), & j=0;\\
y\cdot h(s+j-m-1) , & 1\leq j \leq m;
\end{array}
\right.
\nonumber\\
&&m = m + 1.\}
\end{eqnarray}

\subsubsection{$h(i)\leftarrow\frac{\partial^2 h_{i}(x,y)}{\partial x^2}$ 
used to compute the $x$-component of Laplace term}

\begin{eqnarray}
&&h(j) = 0 ~~ (0\leq j\leq 2, {\rm and} 4\leq j\leq 5);
\nonumber\\
&&h(3) = 2;~~s=3;
\nonumber\\
&&\mathrm{For}~~ m\geq 3 
\nonumber\\
&& \{ s = s + m;
\nonumber\\
&&h(s+j) = \left\{
\begin{array}{cl}
\frac{m}{m-2}\cdot x \cdot h(s-m), & j=0;\\
 y \cdot h(s+j-m-1), & 1\leq j \leq m-2;\\
0, & m-1 \leq j \leq m;
\end{array}
\right.
\nonumber\\
&&m = m + 1.\}
\end{eqnarray}

\subsubsection{$h(i)\leftarrow\frac{\partial^2 h_{i}(x,y)}{\partial y^2}$ used to compute 
the $y$-component of Laplace term}

\begin{eqnarray}
&&h(j) = 0 ~~ (0\leq j\leq 4);
\nonumber\\
&&h(5) = 2;~~s=3;
\nonumber\\
&&\mathrm{For}~~ m\geq 3 
\nonumber\\
&& \{ s = s + m;
\nonumber\\
&&h(s+j) = \left\{
\begin{array}{cl}
0, & 0 \leq j \leq 1;\\
x \cdot h(s+j-m), & 2\leq j \leq m-1;\\
\frac{m}{m-2}\cdot y \cdot h(s-1), & j=m;
\end{array}
\right.
\nonumber\\
&&m = m + 1.\}
\end{eqnarray}

\section{Algorithm test: application to magnetization data}
 
The formula given above are quite general. For application to magnetization data,
it's only needed to replace $X$, $Y$, $Z$ and their corresponding size normalizations 
with $H$, $T$, $M$ and $H=f_{H}(x)$, $T=f_{T}(y)$, $M=f_{M}(z)$. 
Here, we fit the magnetization data of polycrystalline samples 
La$_{1.2}$Sr$_{1.8}$Mn$_{2}$O$_{7}$ obtained with Physical Property Measurement System 
(PPMS) of Quantum Design Company. More details can be found in reference. $^{[4]}$  

\subsection{Fitting without regularization}


\begin{figure}
\includegraphics[ scale=0.5 ]{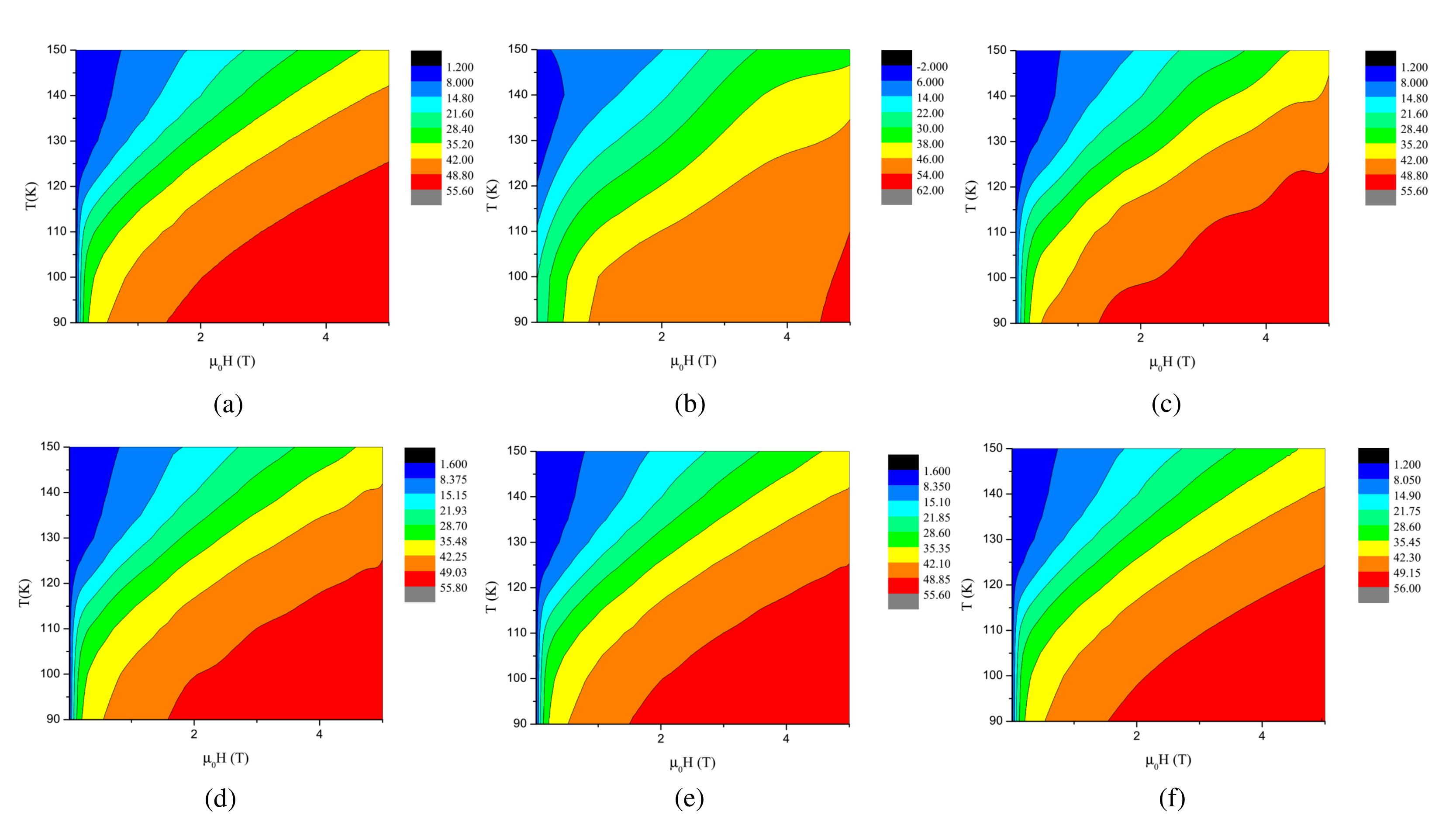}
\caption{
Comparison of experimental data that are processed by the Origin software 
with the fitting effect using two-variable orthogonal polynomials
at regularization parameter $\lambda = 0$.
Shown in (a) are recorded experimental data ($3789$) processed by the Origin software.
By using one half of the whole data ($1895$) as the training group, 
and implementing fitting routines, fitting errors that are achieved and 
the number of orthogonal polynomials that are used list as follow
(b) {9.77E-04}, {15}; 
(c) {9.67E-05}, {49};
(d) {8.75E-06}, {92};
(e) {9.62E-07}, {156};
(f) {1.68E-07}, {350}.
%
}\label{fig1}
\end{figure}


The total number of data used is $3789$. 
With one half of the data ($1895$) uses as the training group, 
the fitting results are shown in Fig. 1.
Note that satisfactory fitting precision can be reached 
using two-variable orthogonal polynomials. 
However, one cannot assess the degree of overfitting 
from the fitting precision, which suggests the necessity to  
introduce cross-validation.
 
The whole data are uniformly sampled according to temperature 
and the sampling factor is set to $3$, namely,
two thirds of the data being put into the training group,
half of the rest one third into the cross-validation group and  
half into the test group.
By setting regularization parameter $\lambda=0$,
we discuss the effect of the number of orthogonal polynomials ($S$) on fitting performance.
The calculated results are shown in Table 1. 
For reference, the overfitting degree corresponding to the test error
is defined in a similar way with that in (\ref{overfit1})  
\begin{equation}
\gamma^{\prime} = \ln \Big | 
\frac{ \sigma_{\rm test}- \sigma_{\rm tr} }{ \sigma_{\rm tr} } 
\Big |.    \label{overfit2}
\end{equation}
It' noted that both $\gamma$ and $\gamma^{\prime}$ reflect fitting performance, 
since the data in the cross-validation and test groups can be interchanged.
With increasing $S$, fitting error decreases and 
overfitting degree increases, which suggests that the overfitting degree 
can be used to monitor fitting performance, 
although it can not always select out the optimal regularization parameter.

\begin{center}
\tabcolsep=15pt  
\small
\renewcommand\arraystretch{1.2}  
\begin{minipage}{15.5cm}{
\small{\bf Table 1.} Fitting errors and overfitting degrees without regularization ($\lambda = 0$).
For comparison, the training error with $S=0$ (i.e., only $h_{0}=1$ is used) is $0.81885{\rm E}-01$.  
 }
\end{minipage}
\vglue5pt
\begin{tabular}{| c | c | c | c | c | c |}  
\hline 
 $S$ & $\sigma_{\rm tr}$ & $\sigma_{\rm cv}$ & $\sigma_{\rm test}$ & $\gamma$ & $\gamma^{\prime}$\\ 
\hline
 {2} & {0.879011E-02} & {0.898435E-02} & {0.887978E-02} & {-3.81} & {-4.58}\\    
 {16} & {0.868127E-03} & {0.752606E-03} & {0.918941E-03} & {-2.02} &{-2.84} \\    
 {50} & {0.966997E-04} & {0.873995E-04} & {0.104586E-03} & {-2.34}  &{-2.51}\\    
 {92} & {0.934213E-05} & {0.158334E-04} & {0.901059E-05} & {-0.36}  &{-3.34}\\   
 {230} & {0.407959E-06} & {0.162667E-04} & {0.624748E-06} & {3.66}  &{-0.63}\\   
\hline
\end{tabular}
\end{center}

\subsection{Fitting with regularization}

\begin{center}
\tabcolsep=15pt  
\small
\renewcommand\arraystretch{1.2}  
\begin{minipage}{15.5cm}{
\small{\bf Table 2.} 
Fitting errors and overfitting degrees at different regularization parameters $\lambda$ 
with the fixed number of orthogonal polynomials $S=78$. 
Parameter $\lambda$ is expressed as $\lambda={\rm e}^{-x}$ for clarity. 
}
\end{minipage}
\vglue5pt
\begin{tabular}{| c | c | c | c | c | c |}  
\hline 
 $x$ & $\sigma_{\rm tr}$ & $\sigma_{\rm cv}$ & $\sigma_{\rm test}$ & $\gamma$ &$\gamma^{\prime}$\\  
\hline
{13} & {0.357571E-02} & {0.346652E-02} & {0.363449E-02} & {-3.49} & {-4.11} \\ 
{15} & {0.193027E-02} & {0.184704E-02} & {0.195522E-02} & {-3.14} & {-4.35}  \\
{17} & {0.903203E-03} & {0.883599E-03} & {0.907491E-03} & {-3.83} & {-5.35}\\
{19} & {0.489599E-03} & {0.535881E-03} & {0.486862E-03} & {-2.36} & {-5.17}\\ 
{21} & {0.200504E-03} & {0.322430E-03} & {0.201699E-03} & {-0.50} & {-5.34}\\
{23} & {0.493752E-04} & {0.845318E-04} & {0.526989E-04} & {-0.34} & {-2.92}\\
\hline
\end{tabular}
\end{center}

Given regularization parameter $\lambda$, 
coefficients of orthogonal polynomials in the fitting expression, $b_{t}$,
are determined according to (\ref{bt2}).
On increasing the number of orthogonal polynomials,
training error $\sigma_{\rm tr}$ as well as cross-validation error $\sigma_{\rm cv}$ 
decreases until terminating the fitting routine according to some criterion. 
Hence, the number of orthogonal polynomials, $S$,  at fixed $\lambda$ 
is automatically identified by the criterion 
that judges where to terminate the fitting routine.

The rising question is how to identify the best $\lambda$?  
Shown in Table 2 are fitting errors 
$\sigma_{\rm tr}$, $\sigma_{\rm cv}$, $\sigma_{\rm test}$), 
and corresponding overfitting degrees $\gamma$, $\gamma^{\prime}$,
at the fixed number of orthogonal polynomials, $S=78$.
The regularization parameter  
$\lambda$ is expressed as $\lambda={\rm e}^{-x}$ with $13\leq x\leq 23$ for clarity.
It's noticed that $\sigma_{\rm tr}$, $\sigma_{\rm cv}$ and $\sigma_{\rm test}$ decrease 
with lowering the stiffness of the surface to fit 
($x$ increases and therefore $\lambda$ decreases). 
Hence the ordinary cross-validation scheme seems not 
quite useful here.
Generally speaking, the overfitting degree ($\gamma$) 
reflects the fitting performance although it does not 
monotonously increase while the stiffness decreases.

It's found that the number of polynomials $S$ affects 
the fitting error not as obviously as $\lambda$ does.
For example, in the case with $S=78$, lowering $\lambda$ 
leads to two orders of change in the fitting error.
With fixed $\lambda$, the fitting error with $S=112$ 
is comparable to that with $S=78$.

Another factor that needs considering is 
the sample factor.
In analysing the dependence of overfitting degree 
on the sample factor, it is noticed that $\sigma_{\rm tr}$,
$\sigma_{\rm cv}$, $\sigma_{\rm tr}$  increase 
with increasing sample factor, namely, 
reducing the size of the training group.  
The overfitting degree decreases 
(with magnitude increases) with the sample factor, suggesting that 
the value of $\sigma_{\rm cv}/\sigma_{\rm tr}$ decreases.
By comparing cases with different sampling factors, 
it's noticed that, in spite of fitting errors as well as 
overfitting degree changing with the sample factor 
the evolving tendencies of the overfitting 
degree with $\lambda$ are similar. 
Thus we can find valuable clues to identify the best $\lambda$.

\begin{center}
\tabcolsep=15pt  
\small
\renewcommand\arraystretch{1.2}  
\begin{minipage}{15.5cm}{
\small{\bf Table 3.} Fitting errors and overfitting degrees at different regularization parameters $\lambda$. 
The number of orthogonal polynomials $S$ is automatically determined by the routine-terminating criterion. 
Parameter $\lambda$ is expressed as $\lambda={\rm e}^{-x}$ for clarity. }
\end{minipage}
\vglue5pt
\begin{tabular}{| c | c | c | c | c | c | c |}  
\hline 
 $x$ & $S$ &$\sigma_{\rm tr}$ & $\sigma_{\rm cv}$ & $\sigma_{\rm test}$ & $\gamma$ &$\gamma^{\prime}$\\  
\hline
{10}& 50 & {0.528493E-02} & {0.510816E-02} & {0.531742E-02} & {-3.40} & {-5.07} \\ 
{20}&78 & {0.329224E-03} & {0.420514E-03} & {0.326656E-03} & {-1.28} & {-4.82}  \\
{30} & 156 & {0.568479E-05} & {0.233510E-04} & {0.524090E-05} & {1.13} & {-2.55}\\
{40} & 201 &{0.696329E-06} & {0.273960E-05} & {0.724544E-06} & {1.07} & {-3.19}\\ 
\hline
\end{tabular}
\end{center}

Table 3 summarizes the optimal fitting precision and corresponding overfitting degrees 
with different regularization parameter $\lambda$. 
The sampling factor is the same to preceding tables.
It's noted that overfitting begins to occur at $\lambda = 30$.
Since the cross-validation and test group can interchange data, 
one needs to compare $\gamma$ with $\gamma^{\prime}$ 
in order to select out the proper $\lambda$.

\subsection{Further Discussions}

Here we discuss aspects that have not been mentioned above. 

First of all, the algorithm in this paper can be further optimized
so that the computing efficiency is increased and memory decreased.
For example, in the re-orthogonalization step of iterating orthogonalization,
the projection of the just orthogonalized polynomial is subtracted 
from all those subsequent polynomials which are not orthogonalized yet;
this definitely increases the the computing cost 
since redundant polynomials are generated to ensure the fitting precision.  

Secondly, if only physical quantities that are 
at experimental-recorded magnetic fields and temperatures      
are concerned, the computing results for the linearly independent functions and 
orthogonal polynomials used to fit can be saved to compute these physical quantities.

Thirdly, after the fitting expression obtained, 
uniform gridding and interpolation of the data can be readily achieved. 
Thus, besides the iterative method used in this work,
operations like  numerical derivatives and integrals can be easily executed.


%
Fourthly, when the coefficients of linear independent functions in 
orthogonal polynomials and those of the orthogonal polynomials in the fitting expression 
are computed, we have not considered the effect of measuring errors.
However, the measuring error is always there. What's more, the measuring error 
of experimental data at extremely weak field is usually bigger. 
Because we have taken into accounted all experimental data with equal weight,
the measuring error under weak field will do harm to global fitting performance.  
To solve this problem, one method is to abandon the data under weak field; 
another method is to introduce a local weight 
factor that depends on magnetic field and temperature so that the 
influence of experimental data is confined within a nearby area.

Fifthly, it's noted that in above fitting, 
rotational symmetry of magnetic systems has not been considered.
Because of this symmetry, the terms involving even-number powers of magnetic field 
should not appear in the analytic expression of magnetization.
While talking about the properties closely associated with this symmetry,
it needs eliminating such terms in the set of linearly independent functions.

Finally, since the algorithm considers only the 
dependence on magnetic field and temperature, 
the effect of other processes such as rotation of crystal grains
and the change of sample volume are not adequate disclosed, that may be 
the reason for larger fitting error at weak fields. 
Another issue not considered here is the influence of demagnetization factor 
whose value in polycrystalline materials is hard computed.

\section{Conclusions} 

To conclude, by adding the Laplace term into the fitting error expression,
the regularization method and corresponding cross-validation scheme 
are introduced to two-variable orthogonal polynomials fitting. 
After applying to magnetization data, it's found that 
the regularization scheme does play its role 
through rapidly suppressing the coefficients of higher-order terms in the fitting expression
and therefore effectively relieving the overfitting problem.
With the aid of the concept of overfitting degree, 
It's also shown that the cross validation scheme can be used to select out 
the proper regularization parameter.  
The influences of sampling parameter and sampling factor are also analysed.  
Thus it offers the quite reliable base 
for the following investigations of the magnetic-entropy-change and 
phase-transition properties of magnetic functional materials.

\section{Acknowledgements}
I would thank Wu Hong-ye for many helps in software usage 
and constructive discussions in developing the method in this work. 
And my thanks also give to Wu Ke-han and Zhou Min for 
providing experimental data before their paper published.

\end{document}